\magnification=\magstep1

\newcount\sec \sec=0
\input Ref.macros
\input math.macros
\forwardreferencetrue
\citationgenerationtrue
\initialeqmacro

\def\ends{{\it Ends}}
\def\vis{vis}


\title{Boundary-connectivity via graph theory}

\author{\'Ad\'am Tim\'ar}
\bigskip

\abstract{We generalize theorems of Kesten and Deuschel-Pisztora about the 
connectedness of the  exterior boundary of a connected subset of $\Z^d$, 
where 
``connectedness" and ``boundary" are understood with respect to various 
graphs on the 
vertices of $\Z^d$. These theorems are widely used in statistical physics and 
related areas of probability.
We provide simple and elementary proofs of their results. 
It turns out that the proper way of viewing these questions is graph 
theory, instead of topology.}



Denote by $\Z^d$ the usual nearest-neighbor lattice on $\Z^d$, i.e., two points 
of $\Z^d$ are adjacent if they differ only in one coordinate, by 1.
Let $\Z^{d*}$ 
be the graph on the same vertex set and 
edges
between every two distinct points that differ in every coordinate by at
most 1.
Say that a set of vertices in $\Z^d$ is
*-connected if it is
connected in the graph $\Z^{d*}$. 

In \ref b.DP/ Deuschel and Pisztora prove that the part of the outer vertex boundary of a finite 
connected subgraph $C$ in $\Z^{d*}$ that is visible from infinity (the {\it exterior boundary})
is *-connected. Earlier, Kesten showed
that 
the set of points in the *-boundary of a connected subgraph 
$C\subset \Z^{d*}$ 
that are $\Z^d$-visible from infinity is connected in $\Z^d$ \ref b.K/. 
Similar results were proved about the case when $C$ is in an $n\times n$ 
box of $\Z^d$ \ref b.DP/, or $\Z^{d*}$ \ref b.H/. See the second 
paragraphs of \ref t.DP/ and \ref t.K/ for the precise statements.

We generalize these results about $\Z^d$ and $\Z^{d*}$ to a very 
general family of pairs of graphs, see \ref l.DP/, \ref t.DP/ and \ref t.K/.
Our method also gives an elementary and short alternative to the 
original proofs for the cubic grid case.
This approach seems to be efficient to treat possible other questions about
the connectedness of boundaries.
Although \ref b.K/ mentions that some use of algebraic topology seems to 
be 
unavoidable, the greater generality (and simplicity) of our proof is a 
result of using purely graph theoretic arguments. Also, it makes slight 
modifications of 
the results (such as considering boundaries in some subset of $\Z^d$ instead of 
boundaries in $\Z^d$) straightforward, while previously one had to go through 
the original proofs and make significant modifications.

In two dimensions, the use of some duality argument makes connectedness of 
boundaries more straightforward to prove.
The lack of duality (that is, the 
correspondance that a cycle in one graph is a separating set in its dual) in higher 
dimensions has been responsible for the increasing difficulty and the role of 
topology. Defining duality in higher dimension led to models such as plaquette 
percolation, where hyperfaces (``duals" of edges) are deleted independently 
with some fixed probability, giving rise to 
random surfaces.

Theorems about connectedness of boundaries have a wide use in probability and statistical physics. To 
list some representatives of the many, connectedness of 
the boundaries in \ref b.DP/ and \ref b.K/ are used in the study of 
Ising, Potts and random cluster models \ref b.Pi/, \ref b.GG/, first 
passage percolation \ref b.K/,
Bernoulli percolation \ref b.KZ/, \ref b.AP/ and random 
walks on percolation clusters \ref b.Pe/, entanglement percolation 
\ref b.GH/, greedy lattice animals \ref b.H/. Understanding connectedness of 
boundaries is an essential part for the use of Peierls estimates, and for 
proving the existence of phase transitions. The fundamental role of these results in many statistical physics arguments makes it important to understand these issues properly.
Our generalizations may help 
extend some of these results to graphs beyond $\Z^d$. This was the case in the simplication of the results of \ref b.BB/ in \ref b.T/, and the main lemma in the latter is the starting point of the current paper (see \ref l.regi/).
Even for the cases where $\Z^d$ is considered, the use of elementary graph theoretic arguments instead of topology adds a lot of flexibility and makes the proofs more accessible.

The graphs we consider can be finite or infinite, but we always assume
that they are locally finite (that is, every vertex has finite degree).
Given a subgraph $H$ of a graph $G$, the {\it inner boundary} of $H$ in $G$ is 
the set of vertices in $H$ that are adjacent to some vertex in $G\setminus H$. 
Similarly, the {\it outer boundary} of $H$ in $G$ is the set of vertices in 
$G\setminus H$ that are adjacent to some vertex in $H$. If $G$ is infinite 
and 
$H$ is finite, the {\it exterior} part
of a boundary (of either type) is the set of vertices in the boundary that are 
starting points of some infinite path with no interior vertex in $H$.
The boundaries we consider
are always taken to be outer boundaries, but our arguments would apply 
just as well for inner boundaries.
By a separating set we always understand a separating set of vertices.
In this paper {\it addition is always understood modulo} 2, and this is 
how we 
define the sums of sets of edges (regarded as vectors over the 2-element 
field). In particular, this defines the generation of cycles 
by other cycles. Let the {\it cycle space} of a graph $G$ be the set of all finite subgraphs such that every vertex has an even degree. It is well known that the cycle space is generated by the set of cycles.

For an arbitrary graph $G$, let $\ends (G)$ be the the set of ends 
in $G$, where an end is an equivalence class of infinite simple paths, two 
being equivalent if they can be connected by infinitely many pairwise disjoint paths. So, $\ends(G)=\emptyset$ iff $G$ is finite, and for 
$G=\Z^d$ we have $|\ends 
(G)|=1$. A path from an end $x$ (or, between an end and a 
vertex $y$) is some path in the equivalence class that defines $x$ (and 
starts
from $y$ respectively). A path between two ends $x,y$ is a biinfinite path $P$ such that for any $v\in P$, $P\setminus v$ consists of a path that belongs to $x$ and a path that belongs to $y$. A {\it separating set} between $x\in V(G)\cup\ends (G)$ and a $y\in V(G)\cup\ends (G)$ 
is a subset of $V(G)$ that every path between $x$ and $y$ intersects. A {\it separating set of edges} between $x\in V(G)\cup\ends (G)$ and a $y\in V(G)\cup\ends (G)$ 
is a subset of $E(G)$ that every path between $x$ and $y$ intersects. An important property of minimal separating sets of edges is that they always split a connected graph into two components (this may not be true for separating sets of vertices).

Given some graph $G$ and a graph $H$ containing $G$, say that a cycle $C$ 
in $G$ is {\it chordal 
in} 
$H$ if any two
points in $C$ are adjacent in $H$. If ${\cal C}$ is a set 
of
cycles in $G$, say that ${\cal C}$ {\it is chordal in} $H$, if every cycle 
in ${\cal C}$ is chordal in $H$. 

The next lemma is the key to our proofs. Similar and slightly weaker versions are in \ref b.BB/ and 
\ref 
b.T/.

\procl l.regi
Let $G$ be some graph, and $\Pi$ a minimal separating  set of edges
between two 
points $x,y\in G\cup\ends(G)$. Let ${\cal C}$ be a set of cycles that 
generate the cycle space of $G$. Then for any partition $(\Pi_1,\Pi_2)$ of $\Pi$, 
there is some cycle $O\in{\cal C}$ that intersects both $\Pi_1$ and $\Pi_2$. 
There is also an $O$ with the above property and such that $|O\cap 
\Pi_2|$ 
is odd.
\endprocl

\proof
If $x$ (or $y$) is an end, define $x'$ ($y'$) to be a vertex such that there is 
a 
path between $x$ and $x'$ ($y$ and $y'$) in $G\setminus \Pi$. Otherwise let 
$x':=x$ ($y':=y$).
Choose paths $P_i$ between $x'$ and $y'$, $i=1,2$, such that $P_i$ does not 
intersect $\Pi_{3-i}$. Such paths exist by the minimality of $\Pi$. There is 
a subset $A\subset {\cal C}$ such that 
$$P_1+P_2=\sum_{C\in A} C.$$
Let $A_1\subset A$ be the set of those cycles that intersect $\Pi_1$, and 
$A_2:=A\setminus A_1$. The previous equation can be written as
$$P_1+\sum_{C\in A_1}C=P_2+\sum_{C\in A_2} C.$$  
The right hand side here does not intersect $\Pi_1$, so it has to intersect 
$\Pi_2$ (since $x'$ and $y'$ 
are the only vertices with odd degree in $P_2+\sum_{C\in A_2} C$, so 
they belong to the 
same component of it). Furthermore, $P_2$ contains 
an odd number of elements from $\Pi_2$, 
and every cycle in $A_2$ contains an even number of elements from $\Pi_2$. 
Thus the total number of elements of $\Pi_2$ in the sum on the right 
side is odd. 
We conclude that the left side (regarded as a subgraph of $G$) has to contain 
some cycle $O$ that 
intersects $\Pi_2$ in 
an 
odd number of edges (since $P_1$ doesn't intersect $\Pi_2$), and 
$O\cap \Pi_1\not=\emptyset$ too, by the definition of $A_1$.\Qed

For a subgraph $C$ of $G$, and $x\in V(G)\cup\ends (G)$, the outer
boundary
of
$C$ visible from $x$ is $\partial_{\vis (x)} (C):=\{y\in V(G):\,y$ is
adjacent to some point in $C$, and there is a path between $x$ and $y$
disjoint from $C\}$. When there are two graphs, $G$ and $G'$ on the same
vertex set, we will also use
$\partial^{G'}_{\vis_{G} (x)} (C):=\{y\in V(G):\,y$ is
$G'$-adjacent to some point in $C$, and there is a $G$-path between $x$
and
$y$
disjoint from $C\}$. 
Hence $\partial^G_{\vis_G (x)} (C)=\partial_{\vis (x)} (C)$. 

Let $B_n$ denote the box induced by $\{1,\ldots,n \}^d$ in 
$\Z^d$.
By a basic 4-cycle of $\Z^d$ we mean the 4-cycle surrounding some 2-face 
in
a unit cube in $\Z^d$. Note that the cycle space of $\Z^d$ has a 
generating set of basic 4-cycles:
think about $\Z^d$ as a Cayley graph
for the free Abelian group. Then the set of basic 4-cycles is the set of
all conjugates of the pairwise commutators of the generating elements,
whose products generate any word equal to the identity --- and
cycles of $\Z^d$ correspond to such words.

The *-connectedness of the $\Z^d$-boundary of a finite connected set in $\Z^{d*}$ is shown 
in \ref b.DP/.  We prove a weaker statement here, assuming that the connected set is from $\Z^d$. We will prove the (generalization of) the original version later in \ref t.DP/, with more assumptions on the underlying graphs.

\procl l.DP
Let $G$ be a graph, and 
$G^+$ be a graph that contains $G$. Suppose that there is a 
generating set ${\Delta}_G$ for the cycle space of $G$ that
is chordal in $G^+$. Then for any connected subset $C$ of $G$ and any $x\in
(V(G)\cup\ends (G))\setminus C$, the set 
$\partial^G_{\vis_G (x)} (C)$ induces a 
connected graph in $G^+$.

In particular, any finite connected subset of $\Z^d$ has a *-connected 
exterior $\Z ^d$-boundary, and if $C\subset B_n$, the outer $\Z^d$-boundary of $C$ 
in any component of 
$B_n\setminus C$ is *-connected.
\endprocl

\proof
Let $\Pi:=\{\{u,v\}\in E(G)\, :\, u\in C, v\in \partial^G_{\vis_G (x)} (C)\}$,
Then $\Pi$ is a minimal separating set of edges in $G$ between $C$ and $x$,
because for every edge $e\in\Pi$ there is a $G\setminus \Pi$-path from $x$ to the endpoint of $e$ in $\partial^G_{\vis_G (x)} (C)$, and appending $e$ to this path we get a path from $C$ to $x$ that intersects $\Pi$ only in $e$. 

Let $\partial^{G}_{\vis_G (x)} (C)=S_1\cup S_2$ be an arbitrary partition. Further, partition $\Pi$ to sets $\Pi _i:=\{\{x,y\}\in E(G)\, :\, x\in C, y\in S_i\}$, $i=1,2$. By \ref l.regi/, there is a cycle $O\in \Delta_G$ such that $O\cap\Pi_1\not =\emptyset$ and $O\cap \Pi_2\not =\emptyset$. Take an edge from each of these intersections, and consider their endpoints in $C$. These are adjacent in $G^+$, since $O$ is chordal, and hence the $G^+$-distance of $S_1$ and $S_2$ is 1. Since the partition to $S_1$ and $S_2$ was arbitrary, we conclude that $\partial^G_{\vis_G (x)} (C)$ is $G^+$-connected.
\Qed

For \ref l.DP/ to hold with a $C$ that is $G^*$-connected but not necessarily $G$-connected (which is the form of the result in \ref b.DP/), we need some extra assumptions on the cycle space. Without those, the conclusion of \ref l.DP/ need not hold, as shown by $G=\Z^2$, $G^+=\Z^{2*}\cup\{\{u,v\}\}$, where $\{u,v\}$ is an edge with endpoints at distance 10 in $\Z^2$, and $C$ we choose to be the $G^+$-connected set induced by the 2-neighborhoods of $x$ and $y$ in $G^+$.

\procl t.DP
Let $G^+$ be a connected graph, and $G$ a connected subgraph of $G^+$. Suppose that there is a 
generating set ${\Delta}_G$ for the cycle space of $G$ that
is chordal in $G^+$, and that for every edge $e\in G^+$ there 
is a cycle $O_e$ in $G^+$ such that $O_e\setminus e\subset G$, and $O_e$ 
is chordal in $G^+$. 
Let $C$ be a 
connected subgraph of 
$G^+$, and $x\in 
(V(G)\cup\ends (G^+))\setminus C$. Then $\partial_{\vis_G(x)}^G (C)$ is connected in $G^+$.

In particular, any finite *-connected subset of $\Z^d$ has a *-connected 
exterior $\Z ^d$-boundary, and if $C\subset B_n$, the outer $\Z^d$-boundary of $C$ 
in any component of 
$B_n\setminus C$ is *-connected.
\endprocl

Note that 
\ref t.DP/ is stronger than the one in \ref b.DP/ even in the $\Z^d$ case: it implies that the 
boundary of 
a 
connected subset of $\Z^d$ is connected in the graph 
$\Z^d\cup\{$edges connecting two points of some basic 4-cycle$\}$, which 
does not follow from the topological proof in \ref b.DP/. This 
strengthening was first shown (for $\Z^d$) in \ref b.GG/).

\proofof t.DP
Define $S:=\partial^{G}_{\vis_G (x)} (C)$. Let $\Pi:=\{\{x,y\}\in E(G)\, :\, x\in C, y\in S \}$,
and $H$ be a graph with $V(H)=V(G)$ and $E(H)=G^+ |_C\cup E(G) $ (here by $G^+ |_C$ we denote the subgraph of $G^+$ induced by $C$).
Then $\Pi$ is a separating set of edges between $C$ and $x$ in $H$, and it is a minimal separating set of edges, because for every edge $e\in\Pi$ there is a path in $G\setminus \Pi$ from $x$ to the endpoint of $e$ in $S$, and appending $e$ to this path we get an $H$-path from $C$ to $x$ that intersects $\Pi$ only in $e$. 

Let $\Delta$ be a generating set for the cycles of $H$, consisting of cycles that are chordal in $G^+$ --- we are going to show the existence of such a $\Delta$. By our assumptions ${\cal H}:=\{O_e\, :\, e\in H\setminus G\}\cup \Delta_G$ consists of cycles that are chordal in $G^+$. 
On the other hand, any cycle $U$ in $H$ is generated by ${\cal H}$, 
because 
$U+\sum_{e\in 
U\setminus G} O_e$ is a 2-regular graph in $G$, and hence it is generated 
by 
${\Delta}_G$.

Let $S=S_1\cup S_2$ be an arbitrary partition. Further, partition $\Pi$ to sets $\Pi _i:=\{\{x,y\}\in E(G)\, :\, x\in C, y\in S_i\}$, $i=1,2$. By \ref l.regi/, there is an $O\in \Delta$ with $O\cap\Pi_1\not =\emptyset$ and $O\cap\Pi_2 \not =\emptyset$. Since $O$ is chordal in $G^+$, we obtain that the $G^+$-distance between  $S_1$ and $S_2$ is 1. Since their choice was arbitrary, $S$ necessarily induces a connected graph in $G^+$. 

The case $G=\Z^d$ follows by choosing ${\Delta}$ to be a 
generating set of  
basic 4-cycles. 
For an edge $e\in\Z^{d*}$, let $O_e$ be a cycle 
such that $O_e\setminus e$
only has edges from a unit cube that contains $e$. 

\Qed

The $\Z^d$ version of the following theorem is due to Kesten. Its proof 
in \ref b.K/ takes a section, 
with references to results from algebraic topology. The similiar statement for the box of $\Z^d$ as $G$ was proved in \ref b.H/ (and it did not follow automatically from Kesten's result).

\procl t.K 
Let $G^+$ be a connected graph, and $G$ a connected subgraph of $G^+$. Suppose that there is a 
generating set ${\Delta}_G$ for the cycle space of $G$ that
is chordal in $G^+$, and that for every edge $e\in G^+$ there 
is a cycle $O_e$ in $G^+$ such that $O_e\setminus e\subset G$, and $O_e$ 
is chordal in $G^+$. 
Let $C$ be a 
connected subgraph of 
$G^+$, and $x\in 
(V(G)\cup\ends (G^+))\setminus C$. Then $\partial^{G^+}_{\vis_G (x)} (C)$ 
is connected 
in 
$G$. 

In particular, if $C\subset \Z^d$ is finite and *-connected, then the 
subset of its exterior outer boundary in $\Z^{d*}$ that is accessible by 
an infinite 
path in $\Z^d\setminus C$ is $\Z^d$-connected. If $C$ is a subset of $ 
B_n$, $x\in B_n\setminus C$, then $\partial^{\Z^{d*}}_{\vis_{\Z^d} (x)} 
(C)$ is 
$\Z^d$-connected.

\endprocl

The first half of the proof is very similar to that of \ref t.DP/. The only difference between the proofs is that we have to define the auxiliary graphs $H$ slightly differently, and that we need some more arguments in \ref t.K/ for the conclusion.

\proofof t.K

Define $S:=\partial^{G^+}_{\vis_G (x)} (C)$. Let $\Pi:=\{\{x,y\}\in E(G^+)\, :\, x\in C, y\in S \}$,
and $H$ be a graph with $V(H)=V(G)$ and $E(H)=G^+ |_C\cup E(G)\cup \Pi  $.
Similarly to the proof of \ref t.DP/, $\Pi$ is a minimal separating set of edges between $C$ and $x$ in $H$, and there exists a
$\Delta$ generating set for the cycles of $H$, consisting of cycles that are chordal in $G^+$.

Let $S=S_1\cup S_2$ be an arbitrary partition. Further, partition $\Pi$ to sets $\Pi _i:=\{\{x,y\}\in E(G)\, :\, x\in C, y\in S_i\}$, $i=1,2$. By \ref l.regi/, there is an $O\in \Delta$ with $O\cap\Pi_1\not =\emptyset$ and $|O\cap\Pi_2 |$ odd.

Suppose first that $O$ contains some vertex $v$ not in $C\cup S$. Let 
$C_{G^+}(x)$ be the component of $x$ in $G^+\setminus \Pi$, and let 
$C_H(x)$ be the component of $x$ in $H\setminus \Pi$. By the chordality of $O$, there is an edge between some vertex $w\in O\cap C$ and $v$. If $v\in C_H(x)$, then this would imply 
$v\in S$, contradicting the assumption on $v$. So suppose $v\not \in C_H(x)$. But the cycle $O_{\{v,w\}}$ is such that every edge of it different from $\{v,w\}$ is in $G$. In particular, there is a $G$-path from $v$ to $S$:
this path goes from $v$ to the element $u$ of $S$ that is the neighbor of $w$ inside $O$, and appending this path to the path from $u$ to $x$ gives that $v$ should be in $C_H(x)$, a contradiction.

Hence $V(O)\subset C \cup S$. Call a set of vertices $B\subset S_2$ in $O$ a {\it block}, if $B$ induces a connected subgraph in $O$ (i.e., a subpath), and it is maximal with this property. Let $I$ be the set of edges in $O$ that have exactly one endpoint in $B$. It is clear by the definition that $O\cap \Pi_2\subset I$, and that $|I|$ is even, since every block contributes two edges to it. 
If there is an edge $e$ in $I$ such that the other endpoint of $e$ is in $S_1$, then the proof is finished: $S_1$ and $S_2$ have distance 1 in $H$ (and hence in $G$, since $H|_S=G|_S$). So suppose not: every $e\in I$ has the form $e=\{x,y\}$ with $x\in S_2$, $y\in V(O)\setminus S\subset C$ (using the assumption $V(O)\subset C \cup S$). That is, $I=O\cap\Pi_2$.  But by the fact that $I$ has an even number of elements, this would contradict the choice of $O$ (that $|O\cap\Pi_2|$ is odd). 

The $\Z^d$ case follows by the same argument as at the end of the proof of \ref t.DP/.

\Qed

\procl r.inner

The proofs of \ref t.K/ and \ref t.DP/ show that the conditions on the cycle spaces of $G$ and $G^+$ 
can be weakened or stated differently: the only important thing is that we can generate the cycle space of $H$ by cycles that are chordal in $G^+$.


\endprocl

\medbreak
\noindent {\bf Acknowledgements.}\enspace
I am very grateful to G\'abor Pete for drawing my attention to the 
subject, and 
for helpful conversations. I also thank Geoffrey Grimmett for his 
comments 
on the manuscript. Finally, I am indebted to an anonymous 
referee for corrections and suggestions, which led to further simplifications of the proofs.
\startbib

\bibitem[AP]{AP} Antal, P. \and Pisztora, \'A. (1996) On the chemical 
distance for supercritical Bernoulli percolation {\it Ann. Probab.} {\bf 
24}, no.2, 1036-1048. 


\bibitem[BB]{BB} Babson, E. \and Benjamini, I. (1999) Cut sets and normed cohomology with applications
to percolation {\it Proc. Amer. Math. Soc.} {\bf 127}, 589-597.

\bibitem[DP]{DP} Deuschel, J. \and Pisztora, \'A. (1996) Surface order large 
deviations for high-density percolation
{\it Prob. Theory and Related Fields} {\bf 104}, 467-482.

\bibitem[GG]{GG} Gielis, G. \and Grimmett, G. (2002) Rigidity of the 
interface in percolation and random-cluster models
{\it
J. Stat. Phys.} {\bf 109}, 1-37. 

\bibitem[H]{H} Hammond, A. (2006) Greedy lattice animals: geometry and 
criticality {\it Ann. Probab.} {\bf 34}, no.2, 593-637.

\bibitem[GH]{GH} Grimmett, G. \and Holroyd, A. (2000) Entanglement in 
percolation {\it Proc. London Math. Soc.} (3) {\bf 81}, No. 2, 485-512.

\bibitem[K]{K} Kesten, H. (1986) Aspects of first-passage percolation, in 
{\it \'Ecole d'\'et\'e de probabilit\'e 
de Saint-Flour XIV}, Lecture Notes in Math 1180, Springer-Verlag, 125-264.

\bibitem[KZ]{KZ} Kesten, H. \and Zhang, Y. (1990) The probability of a 
large finite cluster in supercritical Bernoulli percolation {\it Ann. 
Probab.} {\bf 18}, no.2, 537-555. 

\bibitem[Pe]{Pe} Pete, G. (2008) A note on percolation on $\Z^d$: 
isoperimetric 
profile via exponential cluster repulsion
{\it Elect. Comm. Probab.} {\bf 13}, 377--392. 

\bibitem[Pi]{Pi} Pisztora, \'A.  (1996) Surface order large
deviations for Ising, Potts and percolation models
{\it Probability Theory Rel. Fields} {\bf 104}, 427-466.

\bibitem[T]{T} Tim\'ar, \'A. (2007) Cutsets in infinite graphs {\it 
Combin. Probab. and Comp.} {\bf 16}, issue 1, 159-166.

\endbib
 
\bibfile{\jobname}
\def\noop#1{\relax}
\input \jobname.bbl

\filbreak
\begingroup
\eightpoint\sc
\parindent=0pt\baselineskip=10pt

Hausdorff Center for Mathematics, Universit\"at Bonn, D-53115 Bonn

\emailwww{adam.timar[at]hcm.uni-bonn.de}{}
\htmlref{}{http://www.hausdorff-center.uni-bonn.de/people/timar/}
\endgroup

\bye